\documentclass[12pt]{article}
\usepackage{epsfig,amssymb,epic,eepic,color}

\newcommand{\be}{\begin{enumerate}}
\newcommand{\ee}{\end{enumerate}}
\newcommand{\bi}{\begin{itemize}}
\newcommand{\ei}{\end{itemize}}

\def\R{\mathbb{R}}

\def\be{\beta}

\def\vp{\varphi}

\def\la{\lambda}

\def\si{\sigma}

\def\ep{\varepsilon}

\def\nd{\noindent}

\newenvironment{demo}{{\bf Proof: }}{\hfill $\diamond$ \medskip}

\newtheorem{theo}{Theorem}
\newtheorem{prop}{Proposition}[section]
\newtheorem{stat}[prop]{Statement}

\newtheorem{defi}[prop]{Definition}
\newtheorem{lemm}[prop]{Lemma}

\begin{document}
%\Large
%\huge
\sloppy
\date{\today}
\title{Quasi-energy function for
diffeomorphisms with wild separatrices}
\author{V.~Grines\thanks{N. Novgorod
State University, Gagarina 23, N. Novgorod,
603950 Russia, grines@vmk.unn.ru.}\and F.~
Laudenbach\thanks{Laboratoire de
math\'ematiques Jean Leray,  UMR 6629 du
CNRS, Facult\'e des Sciences et Techniques,
Universit\'e de Nantes, 2, rue de la
Houssini\`ere, F-44322 Nantes cedex 3,
France,
francois.laudenbach@univ-nantes.fr.}\and O.~
Pochinka\thanks{N. Novgorod State
University, Gagarina 23, N. Novgorod, 603950
Russia, olga-pochinka@yandex.ru.}}

\maketitle
\begin{abstract}
According to Pixton \cite{Pi1977}  there are Morse-Smale diffeomorphisms of $\mathbb S^3$ which have no energy function, that is a Lyapunov function whose critical points are all periodic points of the diffeomorphism.   We introduce the concept of quasi-energy function for a Morse-Smale diffeomorphism as a Lyapunov function with the least number of critical points and construct a quasi-energy function for any  diffeomorphism from some class of Morse-Smale diffeomorphisms on $\mathbb S^3$.
\end{abstract}

\nd {\it Mathematics Subject Classification:} 37B25, 37D15, 57M30.

\nd {\it Keywords:} Morse-Smale diffeomorphism, Lyapunov function, Morse theory.

First and third authors thank grant RFBR  No 08-01-00547
of Russian Academy  for partial financial support.

\section{Formulation of results}

According to \cite{GrLaPo}, given a closed smooth $n$-manifold $M^n$ and  a Morse function $\varphi:M^n\to\mathbb R$ is called a {\it Morse-Lyapunov function} for Morse-Smale diffeomorphism  $f:M^n\to M^n$ if:

1) $\varphi(f(x))<\varphi(x)$ if $x\notin Per(f)$ and
$\varphi(f(x))=\varphi(x)$ if $x\in Per(f)$, where $Per(f)$ is the set of periodic points of $f$;

2) any point $p\in Per(f)$ is
a non-degenerate maximum of  $\varphi\vert_{W^u(p)}$ and a
non-degenerate minimum of
$\varphi\vert_{W^s(p)}$.

\begin{defi} Given  a Morse-Smale
diffeomorphism $f:M^n\to M^n$, a   function $\varphi: M^n\to \mathbb R$  is a {\it quasi-energy function} for $f$ if $\varphi$ is a Morse-Lyapunov function for $f$ and has the least possible number of critical points among all Morse-Lyapunov functions for $f$.
\end{defi}

In this paper we consider the  class $G_4$ of
Morse-Smale diffeomorphisms $f:\mathbb S^3\to \mathbb S^3$ whose
nonwandering set consists of exactly four fixed
points: one source $\alpha$, one saddle $\sigma$
and two sinks $\omega_1$ and $\omega_2$.
It follows from \cite{S3} (theorem 2.3), that the
closure of each connected component (separatrix) of
the  one-dimensional manifold $W^u(\sigma)\setminus\sigma$
is homeomorphic to a segment which
consists of this separatrix and two
points:  $\sigma$ and some sink. Denote by  $\ell_1,\ell_2$
the one-dimensional separatrices
containing the respective sinks $\omega_1,\omega_2$ in their closures.
According to \cite{S3}, $\bar\ell_i, i=1,2$  is everywhere smooth except,
 maybe, at
$\omega_i$. So the topological embedding of $\bar\ell_i$
may be  complicated in a neighborhood of the sink.

According to \cite{ArFo}, $\ell_i$ is
called  {\it tame} (or {\it tamely embedded}) if there is a
homeomorphism $\psi_i:W^{s}(\omega_i)\to{\mathbb R}^n$
such that $\psi_i(\omega_i)=O$, where
$O$ is the origin and
$\psi_i(\bar\ell_i\setminus\sigma)$ is a
ray starting from $O$. In the  opposite case $\ell_i$ is called {\it wild}.
It follows from a criterion in \cite{HGP} that the tameness of $\ell_i$ is
equivalent to the existence of a smooth
3-ball $B_{i}$ around $\omega_i$ in any neighborhood of $\omega_i$ such
 that $\ell_i\cap\partial{B}_{i}$
consists of exactly one point. Using lemma 4.1 from \cite{GrLaPo} it
 is possible to make  this criterion more precise in our dynamical setting:
 $\ell_i$ is tame if and only if there is $3$-ball $B_{\omega_i}$
such that
$\omega_i\in f(B_{\omega_i})\subset int~B_{\omega_i}\subset W^s(\omega_i)$
 and $\ell_i\cap\partial{B}_{\omega_i}$
consists of exactly one point.

It was proved in \cite{BoGr2000} that, for every
diffeomorphism $f\in \mathcal{G}_4$, at least one
separatrix ($\ell_1$ say) is tame. It was
also shown that the topological classification of
diffeomorphisms from $\mathcal{G}_4$ is reduced to the
embedding classifications of the separatrix $\ell_2$; hence
there are infinitely many diffeomorphisms from
$\mathcal{G}_4$ which are not topologically conjugate.

To characterize a type of  embedding of $\ell_2$
 we introduce some special Heegaard splitting of $\mathbb S^3$.
Let us recall that
a three-dimensional orientable manifold is %called 
{\it a
handlebody of  genus $g\geq 0 $} if it is obtained from
a 3-ball by an orientation reversing identification of
$g$ pairs of pairwise disjoint 2-discs in  its boundary.
The boundary of such a handlebody is an orientable
surface of genus $g $.

Let $P^+\subset \mathbb S^3$ be a handlebody of genus $g$ such that
 $P^-= \mathbb S^3\setminus int P^+$ is a handlebody (necessarily
of the same genus as $P^+$). Then the pair $(P^+, P^-)$ is a
Heegaard splitting of genus $g$ of $\mathbb S^3$
with Heegaard  surface $S=\partial P^+=\partial P^-$.

\begin{defi} A Heegaard splitting  $(P^+, P^-)$ of \ $\mathbb S^3$ is said to be
adapted to $f\in  G_4$, or $f$-adapted, if:

a) $\overline{W^u(\sigma)}\subset f(P^+) \subset int~P^+$;

b) $W^s(\sigma)$ intersects $\partial P^+$ transversally and
 $W^s(\sigma)\cap P^+$ consists of a unique 2-disc.

An $f$-adapted Heegaard splitting  $\mathbb S^3=P^+\cup P^-$ is said to be  minimal if its genus is minimal among all $f$-adapted splittings.
\end{defi}

\begin{figure}
\begin{center}
\includegraphics[width=0.6\textwidth]{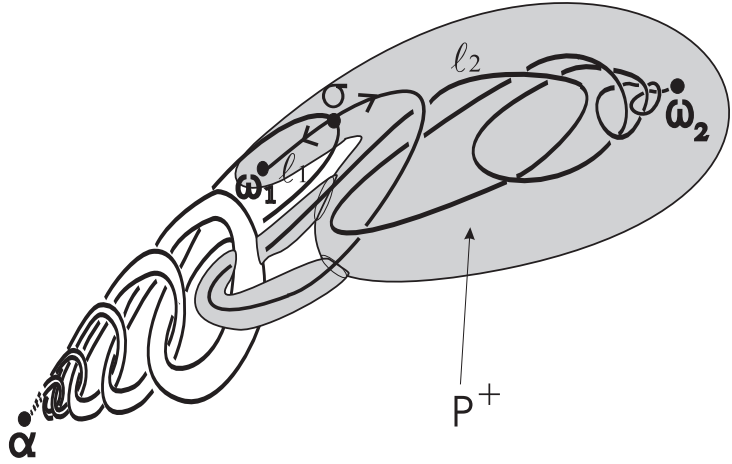}\caption{A diffeomorphism from the class $G_{4,1}$}\label{ld}
\end{center}
\end{figure}

For each integer $k\geq 0$ we denote by $G_{4,k}$ the set of
diffeomorphisms $f\in G_4$ for which the minimal
$f$-adapted Heegaard splitting has genus  $k$.
It is easily seen  that, for  each $f\in G_{4,0}$, $\ell_2$ is tame
and, according to \cite{GrLaPo}, $f$ possesses  an energy
function. Conversely any diffeomorphism in  $G_{4,k},\ k>0,$
has no  energy function (see
\cite{Pi1977}). Figure  \ref{ld} shows the
phase portrait of  a diffeomorphism
 $G_{4,1}$. The main result of this paper is the
following.

\begin{theo} Every quasi-energy function for a diffeomorphism $f\in G_{4,1}$
 has exactly six critical points. \label{qua}
\end{theo}

\section{Recollection of  Morse theory}
\label{MT}

According to Milnor (\cite{Mil65}, section 3),
we use the following definitions.\\

A compact $(n+1)$-dimensional {\it cobordism} is a triad $(W,L_0,L_1)$
where $L_0$ and $L_1$ are closed manifolds of dimension $n$ and
 $W$ is a compact  $(n+1)$-dimensional  manifold whose boundary
consists of the disjoint
union $L_0\cup L_1$. It is an {\it elementary} cobordism
when it possesses a Morse function $\vp: W\to [0,1]$
with only one critical point
and such that $\vp^{-1}(i)=L_i$ for $i=0,1$. When the index  of the unique 
critical point is $r$,
one speaks of an elementary cobordism of index $r$.\\

In this situation, $L_1$ is obtained from $L_0$ by a {\it surgery} of index $r$,
that is: there is an embedding $h:\mathbb S^{r-1}\times \mathbb
D^{n-r+1}\to L_0$
such that $L_1$ is diffeomorphic to the manifold obtained
from $L_0$ by removing the interior of the image of $h$ and gluing
$\mathbb D^r\times \mathbb S^{n-r}$, or
$$L_1\cong\mathbb D^r\times\mathbb S^{n-r}\mathop\bigcup_{h\vert_{\mathbb S^{r-1}\times\mathbb S^{n-r}}}
L_0\setminus int~(h(\mathbb S^{r-1}\times \mathbb D^{n-r+1}))
\,.
$$

Conversely, the following statement holds (see \cite{Mil65}, Theorem 3.12):

\begin{stat} If $L_1$ is obtained from $L_0$ by a surgery of index $r$,
then there exists an elementary cobordism $(W,L_0,L_1)$ of
index $r$.
\label{st4}
\end{stat}

\begin{figure} \begin{center}
\epsfig{file=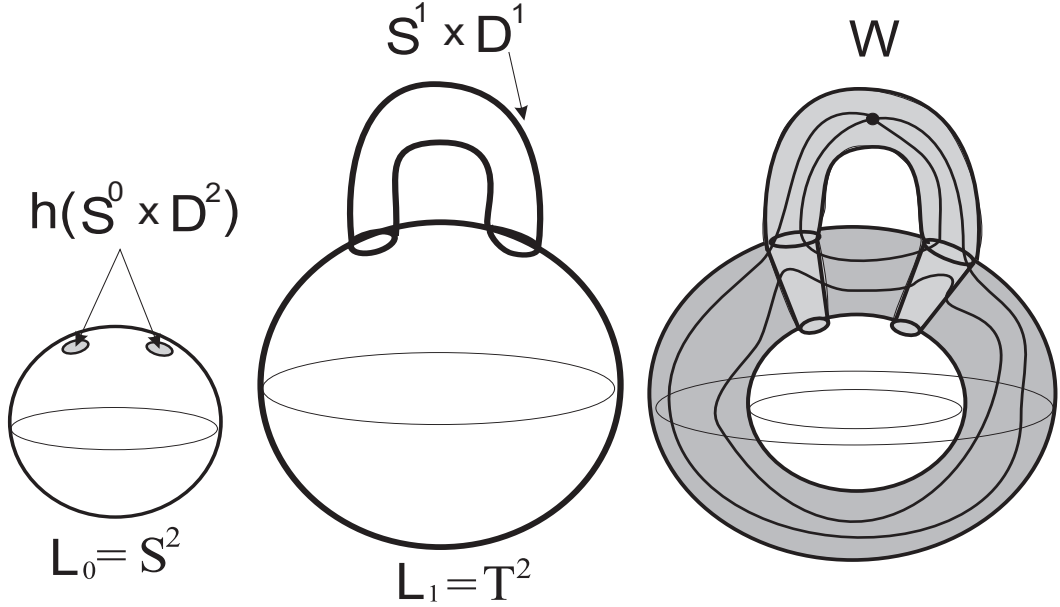, width=14. true cm,
height=8. true cm} \caption{An elementary cobordism}
\label{rec}
\end{center}
\end{figure}

On  figure \ref{rec} it is seen a surgery
of index  $1$ from  the 2-sphere to the 2-torus with some  level sets of a Morse function on the corresponding elementary cobordism.\\

Finally, we recall the weak Morse inequalities
(see \cite{Mil1996}, Theorem 5.2).

\begin{stat} Let $M^n$ be a closed manifold,
$\varphi:M^n\to\mathbb{R}$ be a Morse
function, $C_q$ be the number of
critical points of index $q$ and
$\beta_q(M^n)$ be the $q$-th Betti number
of the manifold $M^n$. Then
$\beta_q(M^n)\leq C_q$ and the Euler characteristic
$\chi(M^n):=\sum\limits_{q=0}^n(-1)^q\beta_q(M^n)$ equals
$\sum\limits_{q=0}^n(-1)^q C_q$.
\label{st10}
\end{stat}

\section{Proof of Theorem \ref{qua}}
Let $f$ be a Morse-Smale diffeomorphism of the 3-sphere  belonging to
 $G_{4,1}$.
As the number of critical points of any Morse function on a closed 3-manifold  is even (it  follows from statement \ref{st10})
and greater than four (as $Per(f)\subset Cr(\varphi)$ and $\ell_2$ is wild)  then, for proving theorem \ref{qua}, it is enough to construct
a Lyapunov function with six critical points.

\subsection{Auxiliary statements}

For the proof of  the following statements  \ref{loc} and \ref{w} we refer
to \cite{GrLaPo}, lemma 2.2 and lemma 4.2.

\begin{stat}  Let $p$ be a fixed point of a Morse-Smale
diffeomorphism $f:M^n\to M^n$ such that
$\dim W^u(p)=q$. Then, in  some neighborhood $U_{p}$ of
$p$, there exist  local coordinates
$x_1,\dots,x_n$ vanishing at $p$ %for which $p=(0,\dots,0)$
and an energy function $\varphi_{p}:U_p\to\mathbb R$
such that
$$\varphi_{p}(x_1,\dots,x_n)=q-x_1^2-\dots-x_q^2+x_{q+1}^2+\dots+x_n^2$$
and $(TW^u(p)\cap U_{p})\subset Ox_1\dots x_q$,
$(TW^s(p)\cap U_{p})\subset Ox_{q+1}\dots x_n$.
\label{loc}
\end{stat}

\begin{stat}  Let $\omega$ be a fixed sink of a
Morse-Smale diffeomorphism $f:M^3\to M^3$ and
$B_\omega$ be a 3-ball with boundary $S_\omega$
such that $\omega\in f(B_\omega)\subset
int~B_\omega\subset W^s(\omega)$. Then there exists
an energy function
${\varphi}_{B_\omega}:B_{\omega}\to\mathbb R$
for $f$ having $S_{\omega}$ as a level set.
\label{w}
\end{stat}

\begin{lemm} Let $\omega$ be a fixed sink of a Morse-Smale diffeomorphism
$f:M^3\to M^3$ and $Q_\omega$ be a solid torus such that
$\omega\in f(Q_\omega)\subset int~Q_\omega\subset W^s(\omega)$.
Then there exists a 3-ball $B_\omega$ such that
$f(Q_{\omega})\subset B_\omega\subset
int~Q_{\omega}$. \label{seq}
\end{lemm}
\begin{demo}  Let $D_0$ be a meridian disk in $Q_\omega$ such
 that $\omega\notin{D_0}$.
As $Q_\omega\subset W^s(\omega)$  there is an integer
 $N$ such that
$f^n(Q_\omega)\cap{D}_0=\emptyset$ for every
$n>N$. We may also
assume that  $D_0$  is transversal to
$G=\bigcup\limits_{n\in\mathbb Z}f^{n}(\partial Q_\omega)$,
 and hence $G\cap int~D_0$ consists of a finite family $\mathcal C_{D_0}$ of
intersection
curves. Each intersection curve $c\in \mathcal C_{D_0}$ belongs
to $f^{k}(\partial Q_\omega)$ for some integer $k\in\{1,\dots,N\}$.
There are two cases: (1) $c$ bounds a disk on $f^{k}(\partial Q_\omega)$;
(2) $c$ does not bound a disk on $f^{k}(\partial Q_\omega)$.
Let us decompose  $\mathcal C_{D_0}$ as union of two pairwise disjoint
 parts $\mathcal C_{D_0}^1$ and  $\mathcal C_{D_0}^2$ consisting of curves
with  property (1) or  (2), accordingly.

Let us show that there is a meridian disk $D_1$ in $Q_\omega$ such that $D_1$
 is transversal to $G$ and $G\cap int~D_1$ consists of
 family $\mathcal C_{D_1}=\mathcal C_{D_0}^2$ of intersection
curves. If $\mathcal C_{D_0}^1=\emptyset$ then ${D_1}={D_0}$. In the
opposite case for any curve $c\in\mathcal C_{D_0}^1$ denote by $d_c$
the disk on $f^{k}(\partial Q_\omega)$ such that $\partial d_c=c$.
Notice that $d_c$ does not contain a curve from the family
$\mathcal C_{D_0}^2$. Then there is
 $c\in \mathcal C_{D_1}$ which is innermost on
$f^{k}(\partial Q_\omega)$ in the sense that the interior of $d_c$ contains
no intersection curves from $\mathcal C_{D_0}$. For  such a curve $c$
denote  $e_c$ the disk on $D_0$ such that $\partial e_c=c$.
As $int~Q_\omega\setminus D_0$ is an open  3-ball then $e_c\cup d_c$ bounds
a unique 3-ball $b_c\subset int~Q_\omega$. Set $D'_c=(D_0\setminus
e_c)\cup d_c$. There is a smooth approximation
$D_c$ of $D'_c$ such that $D_c$ is a meridian disk on $Q_\omega$,
$D_c$  is transversal to $G$. Moreover $G\cap int~D_c$ consists
of a family  $\mathcal C_{D_c}$ of intersection
curves having less elements than $\mathcal C_{D_0}$; indeed, $c$ disappeared
and also all curves from $\mathcal C_{D_0}$ lying in $int~e_c$.
 We will repeat this process until getting a meridian disk $D_1$
 with the required property.

Now let $c\in\mathcal C_{D_1}$, $c\in f^k(\partial Q_\omega)$.
Denote  $e_c$ the disk that $c$ bounds in $D_1$.
 Let us choose $c$ innermost in $D_1$ in the sense that the interior of $e_c$
 contains no intersection curves from $\mathcal C_{D_1}$.
There are two cases: (a) $e_c\subset f^{k}(Q_\omega)$ and
(b) $int~e_c\cap f^{k}(Q_\omega)=\emptyset$.

In case (a) $e_c$ is a meridian disk of $f^{k}(Q_\omega)$ and
$D=f^{-k}(e_c)$ is a meridian
disks in $Q_\omega$  such that $f(Q_\omega)\cap D=\emptyset$. Indeed,
 by construction $int~e_c\cap G=\emptyset$, hence $int~D\cap G=\emptyset$.
Thus we can find the required 3-ball $B_\omega$ inside
$int~Q_\omega\setminus{D}_1$.

In case (b) there is a tubular neighborhood
$V(e_c)\subset int~Q_\omega$ of the disk $e_c$ such that
$G\cap int~V(e_c)=\emptyset$ and  $B_k=f^{k}(Q_\omega)\cup V(e_c)$ is 3-ball.
Then $f^{k}(Q_\omega)\subset B_k\subset int~f^{k-1}(Q_\omega)$.
Thus $B_\omega=f^{1-k}(B_k)$ is the required 3-ball.
\end{demo}

\subsection{Construction of a quasi-energy function
for a diffeomorphism $f\in G_{4,1}$}
As a similar  construction  was done in  section 4.3 of \cite{GrLaPo},   we only  give  a  sketch of it below.

\begin{enumerate}

\item Construct an energy function $\varphi_{p}:U_p\to\mathbb R$
 near each fixed point $p$ of $f$
as in statement \ref{loc}.

\item By definition of the class $G_{4,1}$, for each
$f\in G_{4,1}$ there is a solid torus $P^+$ belonging to a Heegaard splitting
$(P^+,P^-) $ of $\mathbb S^3$ and such
that:

a) $\overline{W^u(\sigma)}\subset f(P^+) \subset
int~P^+$;

b) $W^s(\sigma)$ intersects $\partial P^+$ transversally and
 $W^s(\sigma)\cap P^+$ consists of a unique
2-disk.

\smallskip
\nd As $\mathbb S^3\setminus
\overline{W^s(\sigma)}$ is the disjoint union $W^s(\omega_1)\cup W^s(\omega_2)$, then by  property b),  the
 disk $P^+\cap W^s(\sigma)$ is separating in $P^+$. Moreover there exists
a neighborhood of  $P^+\cap W^s(\sigma)$, such that after  removing it
from $P^+$ we get a 3-ball  $P_{\omega_1}$
and solid torus $P_{\omega_2}$ with the  following properties for each $i=1,2$:

i) $\omega_i\in f(P_{\omega_i})\subset int~P_{\omega_i}\subset W^s(\omega_i)$;

ii) $\partial{P}_{\omega_i}$ is a Heegaard surface
and $\ell_i\cap \partial{P}_{\omega_i}$  consists of exactly one point.\\

\begin{figure}
\begin{center}
\includegraphics[width=0.5\textwidth]{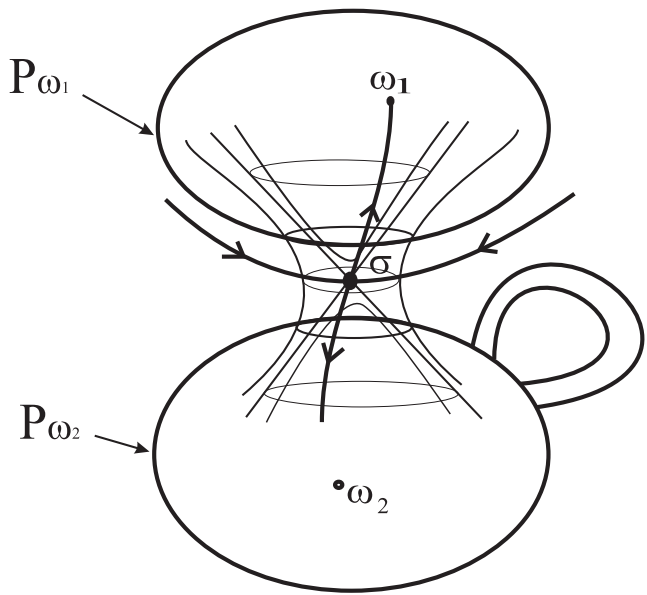}\caption{Heegaard decomposition $(Q^+,Q^-)$}\label{suflu}
\end{center}
\end{figure}

Due to the $\lambda$-lemma\footnote{The $\la$-lemma claims
that $f^{-n}(S_{\omega_i})\cap U_\si$
tends to $\{x_1=0\}\cap U_\si$
in the $C^1$ topology when $n$ goes to $+\infty$.}
(see, for example,\cite{Pa}), replacing $P_{\omega_i}$
 by $f^{-n}(P_{\omega_i})$ for some $n>0$ if necessary,
we may assume that $\partial P_{\omega_i}$ is transversal to the regular part of the
critical level set
$C:=\varphi^{-1}_{\sigma}(1)$ of the function $\varphi_\sigma$
and the intersections  $C\cap\partial P_{\omega_i}$ consist of exactly one circle.
For $\ep\in(0,\frac12)$ define  $H^+_\ep$ as  the closure  of $\{x\in
U_{\si}\mid x\notin (P_{\omega_1}\cup
P_{\omega_2}), \ \varphi_{\sigma}(x)\leq
1+\ep\}$ and set  $P^+_\ep=P_{\omega_1}\cup P_{\omega_2}\cup H^+_\ep$.
In the same way as in  \cite{GrLaPo}
it is possible to choose $\ep>0$ such that $\partial P_{\omega_i}$
intersects transversally   each level set with value in
$[1-\varepsilon$, $1+\varepsilon]$; this intersection consists of  one circle.
Taking a   smoothing $Q^+$ of $P^+_\ep$ we have $f(Q^+)\subset int~Q^+$
and $\Sigma:=\partial Q^+$ is a Heegaard surface of genus 1.
 Let $Q^-$ be the closure of $\mathbb{S}^3\setminus int~Q^+$ (see figure \ref{suflu}). It is easy to check that
 $Q^+$ %, like $P^+_\varepsilon$, (due to the corner, it is not true for it)
  is isotopic to $P^+$. Therefore, the pair  $(Q^+, Q^-)$ is an $f$-adapted Heegaard splitting with the property
that the disk $Q^+\cap W^s(\si)$ lies in $U_\si$.

\item For each $i=1,2$, let $\tilde P_{\omega_i}$  be a handlebody of  genus $i-1$ such that
$f(P_{\omega_i})\subset \tilde P_{\omega_i}\subset int~P_{\omega_i}$,
$\partial\tilde P_{\omega_i}$ intersects transversally each level set
 with value in $[1-\varepsilon,1+\varepsilon]$ along
 one circle and $P_{\omega_i}\setminus int~\tilde P_{\omega_i}$ is
diffeomorphic to $\partial P_{\omega_i}\times[0,1]$.
 Define $d_i$ as the closure of
$\{x\in U_{\sigma}\mid x\in (W^s(\omega_i)\setminus\tilde P_{\omega_i}),
\ \vp_{\sigma}(x)=1-\ep\}$. By  construction $d_i$ is a disk whose boundary
curve bounds a disk $D_i$ in $\partial \tilde P_{\omega_i}$.
 We form $S_i$ by removing the interior of  $D_i$ from
$\partial \tilde P_{\omega_i}$  and gluing the $d_i$.
 Denote $P(S_i)$ the handlebody of genus $i-1$ bounded by $S_i$ and containing
$\omega_i$. As in    \cite{GrLaPo}  it is possible
to choose $\varepsilon$ such that $f(P(S_i))\subset int~P(S_i)$.

Let $K$ be the domain between $\partial Q^+$ and $S_1\cup S_2$.
We introduce $T^+$, the closure  of
$\{x\in \mathbb S^3\mid x\notin (P_{\omega_1}\cup P_{\omega_2}),\ 1-\ep\leq\varphi_{\sigma}(x)\leq
1+\ep\}$; observe $T^+\subset U_\si$. We define a function $\vp_{_K}:
K\to \R$ whose value is $1+\ep$ on $\partial Q^+$, $1-\ep$ on $S_1\cup S_2$,
coinciding with  $\vp_{\sigma}$ on $K\cap T^+$ and without critical points
outside $T^+$. This last condition is easy to satisfy as the domain
in question is a product cobordism. In a similar way to \cite{GrLaPo}, section 4.3,
one can check that  %$\vp^+$ This function is not defined
$\vp_{_K}$ is a Morse-Lyapunov function.

\item As $P(S_1)$ is a 3-ball such that
$\omega_1\in f(P(S_1))\subset
int~P(S_1)\subset W^s(\omega_1)$, then by statement \ref{w}
there is an energy function
$\varphi_{_{P(S_1)}}:P(S_1)\to\R$
for $f$ with $S_1$ as a level set with value
$1-\varepsilon$.

\item As $P(S_2)$ is a solid torus such that
 $\omega_2\in f(P(S_2))\subset
int~P(S_2)\subset W^s(\omega_2)$,
then according to lemma \ref{seq} there is a 3-ball $B_{\omega_2}$ such that
$f(P(S_2))\subset B_{\omega_2}\subset
int~P(S_2)$. As in  the previous item, there is  an energy
function $\varphi_{_{B_{\omega_2}}}:B_{\omega_2}\to\R$ for $f$ with
$\partial{B_{\omega_2}}$ as a level set with value $\frac12$.

\item As $P(S_2)$ is a solid torus, it   is obtained from
a 3-ball by an orientation reversing identification of a pair of disjoint 2-discs in  its boundary; hence the solid torus  is the union of a 3-ball and an elementary cobordism of index 1. Since,  up to isotopy, there is only one 3-ball in the interior of a solid torus,  then $(W_{\omega_2},\partial B_{\omega_2},S_2)$ is an
elementary cobordism
 of
index $1$, where 
 $W_{\omega_2}=P(S_2)\setminus int~B_{\omega_2}$. Hence
  $W_{\omega_2}$ possesses a Morse function
$\vp_{_{W_{\omega_2}}}$ with only one critical point of index 1 and such that
 $\vp_{_{W_{\omega_2}}}(\partial B_{\omega_2})=\frac12$,
$\vp_{_{W_{\omega_2}}}(S_2)=1-\varepsilon$.

\item Define the smooth function ${\varphi}^+:Q^+\to
\mathbb{R}$ by the formula

${\varphi}^+(x)=\cases{\varphi_{_K}(x),
x\in K;
\cr\varphi_{_{P(S_1)}}(x),x\in P(S_1);
\cr\varphi_{_{B_{\omega_2}}}(x),x\in B_{\omega_2};
\cr\varphi_{_{W_{\omega_2}}}(x),x\in W_{\omega_2}.\cr}$

Then ${\varphi}^+$ is a Morse-Lyapunov function for $f\vert_{Q^+}$ with one
additional critical point.

\item By the construction $Q^-$ is a solid torus such that
 $\alpha\in f^{-1}(Q^-)\subset int~Q^-\subset W^u(\alpha)$. Since $\alpha$ is a
sink for $f^{-1}$ then,  as in item 4, there is a 3-ball $B_{\alpha}$ such that
 $f^{-1}(Q^-)\subset B_{\alpha}\subset
int~Q^-$ and an energy function
$\varphi_{_{B_{\alpha}}}:B_{\alpha}\to\R$ for $f^{-1}$ with
$\partial{B_{\alpha}}$ as a level set of value $\frac12$.

\item Similarly to item 5, $\partial Q^-$ is obtained from $\partial B_{\alpha}$
by a surgery of index $1$.
Therefore $(W_\alpha,\partial Q^-,\partial B_\alpha)$ is an elementary
cobordism  of
index $1$, where $W_\alpha=Q^-\setminus int~B_{\alpha}$. Hence, $W_\alpha$ possesses a Morse
function $\vp_{_{W_\alpha}}$ with only one critical point of index 1. We may choose
 $\vp_{_{W_\alpha}}(\partial B_{\alpha})=\frac12$,
$\vp_{_{W_\alpha}}(\partial Q^-)=2-\varepsilon$.

\item Define the smooth function ${\varphi}^-:Q^-\to
\mathbb{R}$ by the formula

${\varphi}^-(x)=\cases{3-\varphi_{_{B_{\alpha}}}(x),
x\in \varphi_{_{B_{\alpha}}};
\cr 3-\vp_{_{W_\alpha}}(x),x\in \vp_{_{W_\alpha}}.\cr}$

Then ${\varphi}^-$ is a Morse-Lyapunov function for $f\vert_{Q^-}$
 with one additional critical point.

\item The function $\varphi:\mathbb S^3\to \mathbb R$ defined by
$\varphi\vert_{Q^+}=\varphi^+$ and
$\varphi\vert_{Q^-}=\varphi^-$ is the required
Morse-Lyapunov function for  the diffeomorphism $f$
with exactly six critical points.

\end{enumerate}

\section*{\addcontentsline{toc}{section}{Reference}}

\end{document}